\documentclass{amsart}%[12pt] before {amsart}
\newtheorem{theorem}{Theorem}%[section]
\newtheorem{proposition}{Proposition}%[section]
\newtheorem{lemma}{Lemma}

\theoremstyle{definition}

\theoremstyle{remark}

\numberwithin{equation}{section}

\newcommand{\scr}{\scriptstyle}
\def\sumprime_#1{\setbox0=\hbox{$\scriptstyle{#1}$}
\setbox2=\hbox{$\displaystyle{\sum}$}
\setbox4=\hbox{${}'\mathsurround=0pt$}
\dimen0=.5\wd0 \advance\dimen0 by-.5\wd2
\ifdim\dimen0>0pt
\ifdim\dimen0>\wd4 \kern\wd4 \else\kern\dimen0\fi\fi
\mathop{{\sum}'}_{\kern-\wd4 #1}}

\begin{document}
\title{Notes on Pair Correlation of Zeros and Prime Numbers}
\author{D. A. Goldston}
\address{Department of Mathematics, San Jose
State University, San Jose, CA 95192, USA}
\email{goldston@math.sjsu.edu}
\thanks{The  author was supported by the Focused Research Group grant (0244660) from NSF and the American Institute of Mathematics. This paper will appear in \lq\lq Recent Perspectives in Random Matrix Theory and Number Theory" to be published by Cambridge University Press.}

%\subjclass{Primary 11N05 ; Secondary 11P32}

\date{\today}

%\keywords{prime number}

%\begin{abstract} 
% \end{abstract}

\maketitle
%\markright{Newton Lecture Notes}
{\footnotesize These notes are based on my four lectures given at the Newton Institute in April 2004 during the Recent Perspectives in Random Matrix Theory and Number Theory Workshop. Their purpose is to introduce the reader to the analytic number theory necessary to understand Montgomery's work on the pair correlation of the zeros of the Riemann zeta-function and subsequent work on how this relates to prime numbers. A very brief introduction to Selberg's work on the moments of $S(T)$ is also given.  }

\section{Introduction and Some Personal History}
In 1973 Montgomery's paper \cite{Mo}, \lq \lq The Pair Correlation of Zeros of the Zeta Function" appeared in the AMS series of Proceedings of Symposia in Pure Mathematics, and a new field of study was born --- slowly. I first came across this paper in 1977, and was probably the only person at Berkeley to read it. Most zeta-function people (as some of us refer to ourselves) recognized the importance of this work and the new phenomena discovered, but it was not clear what to do next. At first, the main interest was in using Montgomery's conjectures to refine the classical results on primes obtained assuming the Riemann Hypothesis. Gallagher and Mueller \cite{GaMu} wrote an important paper on this in 1978, followed by further results from Heath-Brown \cite{HB}. In 1981 I wrote my Ph.D.  thesis on this topic.  A few years later Montgomery and I \cite{GoMo} obtained an equivalence between the pair correlation conjecture and primes.  However, this work attracted little attention --- probably because the results were obtained using Montgomery's conjectures. Then, in the early 1980's everything changed: Odlyzko \cite{O} computed statistics on the zeros and convinced even the most skeptical that after almost a century of intensive study a totally new, unsuspected, and fundamental property of the zeta-function had been discovered. The field has since had a flood of activity, with the generalization of Montgomery's work to higher correlations by Hejhal \cite{He} and Rudnick-Sarnak \cite{RS}, the interpretation of these results in terms of  mathematical physics by Berry and Bogomolny-Keating \cite{BKI,BK2}, the function field case of Katz-Sarnak \cite{KS}, the random matrix model for moments of the zeta-function of Keating and Snaith culminating in \cite{CFKRS}, and a profusion of new work.

In these notes,  I will discuss  Montgomery's results and their relations to primes.  As a unifying tool,  I will use Montgomery's explicit formula  \cite{Mo} to prove a number of  later results that were originally obtained by other methods. This approach was first made use of in part of my Ph.D. thesis, and was based on a suggestion of Montgomery in a letter. At that time Heath-Brown had just finished his paper which covered the same ground, and  I saw no need to publish this material beyond the summary that appeared in \cite{Gtex}. My goal, in line with the emphasis of the workshop on reaching out to beginners in the field, is to provide some of the main ideas used without technicalities and at the same time supply simple details which would be accepted without comment by experts. I have intentionally left out many things to keep these notes focused. The last section on Selberg's theory of $S(T)$ and $\log \zeta(s)$ is somewhat different from the previous ones, and I have decided to state only the main results and present a few of the ideas that are used.

I would like to thank Andrew Ledoan for the many improvements he suggested for these notes.

\section{Basic Facts and Notation}
Following Riemann, we use the complex variable $s= \sigma +i t$.
The Riemann zeta-function $\zeta(s)$ is defined, for $\sigma >1$, by either the Dirichlet series or the Euler product
\begin{equation} \zeta(s) = \sum_{n=1}^\infty \frac{1}{n^s} = \prod_p\left(1-\frac{1}{p^s}\right)^{-1}.\label{2.1} \end{equation}
Here $p$ will always denote a prime, so the product is over all the prime numbers.
To extract information about primes from the Euler product, we compute the logarithmic derivative of the zeta-function and use the power series for $-\log(1-z)$, to obtain, for $\sigma >1$,
\begin{equation} \frac{\zeta'}{\zeta}(s) := \frac{\zeta'(s)}{\zeta(s)} =\frac{d}{ds}\log \zeta(s) = \frac{d}{ds}\Big(\sum_{m=1}^\infty \sum_p\frac{1}{mp^{ms}} \Big) =  -\sum_{n=1}^\infty \frac{\Lambda(n)}{n^s}, \label{2.2}\end{equation}
where the von Mangoldt function $\Lambda(n)$ is given by
\begin{equation}
	\Lambda(n) = 
		\left\{ 
		\begin{array}{ll}
       	\log p, &\mbox{if $n= p^m$, $p$ prime, $m\ge 1$},\\
                  0, & \mbox{otherwise}.
		\end{array}
		\right. 
	\label{2.3} \end{equation}
 The Chebyshev function $\psi(x)$ is the counting function for $\Lambda(n)$ given by
\begin{equation} \psi(x) = \sum_{n\le x}\Lambda(n) .\label{2.4}\end{equation}
Because of the simple relationship with the zeta-function, it is preferable to use $\Lambda(n)$ in place of the indicator function for the primes, and $\psi(x)$ in place of the counting function $\pi(x)$ for the number of primes up to $x$. If needed, one can usually recover $\pi(x)$ from $\psi(x)$  by simple arguments. 
The Prime Number Theorem (PNT) states that  as $x\to \infty$  
\begin{equation} \psi(x) \sim x, \quad \text{or} \quad \pi(x) \sim \frac{x}{\log x}.\label{2.5}\end{equation}
The PNT with the error term obtained by de la Vall\'ee Poussin in 1899 is, for  a small constant $c$,
\begin{equation} \psi(x) = x + O\left(xe^{-c\sqrt{\log  x}}\right) ,\label{2.6} \end{equation}
which on returning to $\pi(x)$ gives (c may differ from equation to equation)
\begin{equation} \pi(x) = \mathrm{li}(x) + O\left(xe^{-c\sqrt{\log  x}}\right) ,\label{2.7} \end{equation}
where the logarithmic integral
\[ \mathrm{li}(x) = \int_2^x \frac{du}{\log u} \]
is the actual main term in the theorem. For the error term above, we have for any constant $A>0$ 
\begin{equation} e^{-c\sqrt{\log  x}} \ll \frac{1}{(\log x)^A}.\label{2.8}\end{equation}
Here,  the Vinogradov notation $\ll$ is equivalent to \lq \lq big oh" of the right-hand side.  This  estimate is freely used when the PNT is invoked. 

We frequently need the Dirichlet series for $\zeta(s)^{-1}$, which
from the Euler product is 
\begin{equation} \frac{1}{\zeta(s)} = \sum_{n=1}^\infty \frac{\mu(n)}{n^s}\label{2.9}\end{equation}
for $\sigma >1$, where the
M\"obius function is defined by $\mu(1)=1$  and
\begin{equation} 
	\mu(n) = 
		\left\{ 
		\begin{array}{ll}
       	(-1)^m, &\mbox{if $n= p_1p_2\cdots p_m$, $p_i$'s distinct},\\
                  0, & \mbox{if $p^2|n$, some $p$}.
		\end{array}
		\right. 
	\label{2.10} \end{equation}

The zeta-function has a simple pole with residue 1 at $s=1$, trivial zeros at $s=-2n$, $n=1,2,3,\ldots$, and  complex zeros 
\begin{equation} \rho = \beta +i \gamma, \quad 0 < \beta <1. \label{2.11}\end{equation}
The inequality $\beta <1$ is the key result needed in the analytic proofs of the PNT. The zeros are positioned symmetrically with the real line and the \lq \lq  $\frac{1}{2}$-line" $\frac{1}{2} +it$, so that $\rho$, $\bar{\rho}$, $1-\rho$, and $1-\bar{\rho}$ are all zeros. The Riemann Hypothesis (RH) is the conjecture that $ \beta = \frac{1}{2}$, and thus $\rho = \frac{1}{2} +i\gamma$. For example, the first 6 zeros in the upper half of the critical strip are
\begin{equation}\begin{split} &\frac{1}{2}+i 14.13472\ldots,\ \ \frac{1}{2}+i 21.02203\ldots,\ \ \frac{1}{2}+i 25.01085\ldots,\\ & \frac{1}{2}+i  30.42487\ldots,\ \ \frac{1}{2}+i 32.93506\ldots,\ \ \frac{1}{2}+i 37.58617\ldots .\label{2.12}\end{split}\end{equation}

To count the number of complex zeros in a given region,  we define 
\begin{equation} n(T) = \left| \{ \gamma :  0<\gamma \le T \} \right|,\qquad  N(T) = \frac{n(T+0)+n(T-0)}{2},\label{2.13} \end{equation}
where $|A|$ denotes the number of elements of the set $A$. Note that $N(T)$ counts any zeros $\gamma =T$ with weight one-half, which arises naturally in the theory;  therefore we always use $N(T)$ in preference to $n(T)$.   The Riemann-von Mangoldt formula for $N(T)$, obtained by applying the argument principle to $\zeta$ and using the functional equation (see \cite{Da}, \cite{In}, \cite{T}), is
\begin{equation}  N(T) = \frac{T}{2\pi}\log \frac{T}{2\pi e} + \frac{7}{8} +R(T) +S(T), \label{2.14}\end{equation}
where 
\begin{equation} R(T) \ll  \frac{1}{T}\label{2.15} \end{equation}
and
\begin{equation} S(T) = \frac{1}{\pi} \arg \zeta\left(\frac{1}{2} +iT\right) \ll \log T.\label{2.16}\end{equation}
In fact, $R(T)$  is continuous, differentiable,  and can be expanded into a series in inverse powers of $T$.
We see that \eqref{2.14} provides a remarkably precise formula for the number of zeros up to height $T$, with the finer details of the vertical distribution of zeros wrapped up in the study of $S(T)$.  In particular we have
\begin{equation} N(T) \sim \frac{T}{2\pi}\log T. \label{2.17} \end{equation}
Another consequence of \eqref{2.14} -- \eqref{2.16} which we make frequent use of is the sharp estimate
\begin{equation} N(T+1)-N(T) = \sum_{T<\gamma \le T+1}1  \ll \log T .\label{2.18} \end{equation}

\section{Explicit Formulas}
To study the relationship between the zeros of the zeta-function and primes you need to be able to work with explicit formulas. There are many such formulas but the best known is the Riemann-von Mangoldt explicit formula, which states that, for $x>1$, 
\begin{equation} \psi_0(x) = x - \sum_{\rho} \frac{x^\rho}{\rho} - \log 2\pi - \frac{1}{2}\log\left(1-\frac{1}{x^2}\right),\label{3.1}\end{equation}
where $\psi_0(x) = \frac{1}{2}(\psi(x+0)+\psi(x-0))$. By \eqref{2.14}, the sum is not absolutely convergent and the terms are added with $\rho$ and $\bar{\rho}$ grouped together.  The explicit formula also contains this information, since on taking $x\to 1^+$ and letting $x=e^u$ we see that
\begin{equation} \sum_{\rho} \frac{e^{\rho u}}{\rho}  = \frac{1}{2}\log \frac{1}{u} +O(1) \quad \textrm{as} \ \ u\to 0^+. \label{3.2}\end{equation}
For applications we usually use the truncated version of \eqref{3.1},
\begin{equation} \psi(x) = x - \sum_{|\gamma|\le T} \frac{x^\rho}{\rho}  +O\left( \frac{x}{T}(\log xT)^2\right)+O\left((\log x) \min\Big(1,\frac{x}{T ||x||}\Big)\right),\label{3.3}\end{equation}
where $||x||$ denotes the distance of $x$ to the closest integer, and the last term reflects the jumps of $\psi(x)$ at the primes and prime powers.
As an example of an application of \eqref{3.3}, assuming  RH we have 
\[ \frac{x^\rho}{\rho} \ll \frac{x^{\frac{1}{2}}}{|\gamma|},\]
and so by \eqref{2.12} and \eqref{2.18}, with $[x]$ denoting the integer part of $x$, 
\[ \sum_{|\gamma|\le T} \frac{1}{|\gamma|} =2 \sum_{1<\gamma\le T} \frac{1}{\gamma}
 \le \sum_{n=1}^{[T]+1}\frac{1}{n}\sum_{n<\gamma \le n+1}1 
 \ll \sum_{n\le 2T}\frac{\log 2n}{n} 
\ll (\log T)^2.\]
Thus, taking $T=x$ in \eqref{3.3} we have 
\begin{equation} \psi(x)= x +O\left(x^{\frac{1}{2}} (\log x)^2\right) .\label{3.4}\end{equation}
It can also be proved that this estimate implies RH, and therefore is equivalent to the RH.\footnote{Even the seemingly weaker estimate $\psi(x)=x+O\left(x^{\frac{1}{2}+\epsilon}\right)$ for any $\epsilon>0$ is equivalent to RH.} Equation \eqref{3.4} is due to von Koch in 1901 and has never been improved. 

We next apply \eqref{3.4} to the problem of large gaps between primes. Let $p_n$ denote the $n$-th prime number. The highest power of a prime $\le x$ is the largest $k$ for which $2^k\le x$, so $k = [\log_2x]$. By the PNT,
\[\begin{split} \psi(x) &= \sum_{p\le x}\log p  + \sum_{2\le m \le\log_2 x} \sum_{p^m\le x} \log p \\ & = \sum_{p\le x}\log p + \sum_{p\le \sqrt{x}}\log p + O\left(\pi\left(x^{\frac{1}{3}}\right)(\log x)^2\right)\\ &
                         = \sum_{p\le x}\log p + O\left(x^{\frac{1}{2}}\right). \end{split}\]
Thus \eqref{3.4} continues to hold when we only sum over primes. For $1\le h \le x$, we have by \eqref{3.4} and differencing that
\[ \sum_{x<p\le x+h}\log p = h +O\left(x^{\frac{1}{2}}(\log x)^2\right).\]
On taking $h = C x^{\frac{1}{2}}(\log x)^2$, with the constant $C$ being larger than the implicit absolute constant in the error term, we conclude that the sum on the left is positive and $\gg h$, and thus the interval $(x,x+h]$ must contain $\gg \frac{h}{\log x}$ primes. If $p_n$ is the first prime in $(x,x+h)$, then 
\begin{equation}  p_{n+1}-p_n < h  \ll {p_n}^{\frac{1}{2}}(\log p_n)^2 .\label{3.5}\end{equation} 

An explicit formula that also exhibits the close connection between zeros and primes is the Landau formula, which states that (for $x$ fixed)  as $T\to \infty$
\begin{equation} \sum_{0<\gamma \le T}x^{\rho} =  -\frac{T\Lambda(x)}{2\pi} +O(\log T) .\label{3.6}\end{equation}
Here we define $\Lambda(x)$ to be zero for real non-integer $x$. Formally this is obtained by differentiating \eqref{3.1} with respect to $x$. The exponential sum over the zeros encodes the information on which integers are primes or prime powers.  Equation \eqref{3.6} is not particularly useful, but Fujii \cite{Fu} and independently Gonek \cite{Gonek} have developed uniform versions which can be used in applications.

An explicit formula of at least historic interest is the Cram\'er explicit formula, which states that for Im$(z)>0$ 
\begin{equation} \begin{split} \sum_{\gamma>0}e^{\rho z}= \frac{e^z}{2\pi i} &\sum_{n=2}^\infty \frac{\Lambda(n)}{n}\left(\frac{1}{z-\log n} +\frac{1}{\log n}\right) \\&
+ \frac{1}{2\pi i} \sum_{n=2}^\infty \frac{\Lambda(n)}{n}\left(\frac{1}{z+\log n} -\frac{1}{\log n}\right) \\ &
+\left(\frac{1}{4}+\frac{ \gamma +\log 2\pi }{2\pi i}\right)\left(1 +\frac{1}{z}\right) + \frac{1}{2\pi i}\frac{\Gamma '}{\Gamma}\left(\frac{z}{2\pi i}\right) + \frac{1}{2}e^z\\& - \frac{z}{2\pi i}\int_0^1e^{sz}\log |\zeta(z)|\, ds -\frac{1}{2\pi iz}\int_0^\infty \frac{t}{e^t-1}\frac{dt}{t+z}.\label{3.7}
\end{split}\end{equation}
On taking
 $z=-\log \tau+iy$, $0<y\le 1$, and letting $\tau \to \infty$ Cram\'er \cite{Cr} proved that 
\begin{equation}-2\pi \text{Re}\sum_{\gamma>0}e^{\rho(-\log \tau +iy)} = \sum_{n=2}^\infty\frac{\Lambda(n)}{n}\left(\frac{y}{\left(\log \frac{\tau}{n}\right)^2+y^2}\right) - \pi +O\left(\frac{1}{\log \tau}\right).\label{3.8} \end{equation}
He used \eqref{3.8} and related formulas in a series of papers starting in 1920 to prove results on primes. One such result is that on RH
\begin{equation} p_{n+1}-p_n \ll {p_n}^{\frac{1}{2}}\log p_n .\label{3.9}\end{equation}
This only saves a logarithm over the trivial use of \eqref{3.4} in \eqref{3.5}  but is the best result known on RH. We will later assume a much stronger hypothesis and only improve \eqref{3.9} by a half-power of a logarithm. On the other hand, Cram\'er conjectured \cite{Cr} that the gaps between consecutive primes are always much smaller than this size. Recent work indicates that Cram\'er's original conjecture may be slightly too strong, but all evidence still suggests
\begin{equation} p_{n+1}-p_n \ll (\log p_n)^2 .\label{3.10}\end{equation}
At one time I had a fondness for Cram\'er's formula and made use of it in my thesis, but I later decided that nothing was to be gained by its use except complicated arguments. The proof of \eqref{3.9}, for instance, can now be done from a smoothed version of  \eqref{3.1} in just a few lines. However, there have been a number of recent papers on the structure of Cram\'er's formula (see \cite{Il}).

Most of these explicit formulas are based on evaluating the contour integral
\[ \mathcal{I}= \frac{1}{2\pi i}\int_{c-i\infty}^{c+i\infty} -\frac{\zeta '}{\zeta} (s) K_z(s)\, ds ,\]
where the kernel $K_z(s)$ is a meromorphic function.  Frequently $K_z(s) =K(s+z)$ or $K_z(s) =K(zs)$. If $c>1$  the Dirichlet series for $\frac{\zeta '}{\zeta}(s)$ converges absolutely and
\[ \mathcal{I} = \sum_{n=2}^\infty \Lambda(n)\hat{K}_z(n), \quad  \hat{K}_z(n) =  \frac{1}{2\pi i}\int_{c-i\infty}^{c+i\infty} K_z(s)n^{-s}\, ds .\] 
One then obtains an explicit formula by moving the contour to the left, thus encountering poles at $s=1$ and at the zeros of $\zeta(s)$, as well as any poles of $K_z(s)$.

Another explicit formula frequently used is the Weil explicit formula,
which contains a general weight function and has the advantage of allowing the relationships between terms to be explicitly exhibited. Our approach here, however, is to use explicit formulas only as tools for studying zeros and primes. Therefore we will take the opposite path and stay specific. 
The formula we will base our work on is due to 
Montgomery  \cite{Mo}. 
\begin{proposition}  Assume the Riemann Hypothesis. For $x\ge 1$,
\begin{equation} \begin{split} 2 x^{\frac{1}{2}-it}\sum_\gamma \frac{x^{i\gamma }}{1+(t-\gamma)^2} &= - \sum_{n=1}^\infty \frac{\Lambda(n) a_n(x)}{n^{it}} +\frac{2 x^{1-it}}{\left(\frac{1}{2}+it\right)\left(\frac{3}{2}-it\right)}\\& +x^{-\frac{1}{2}}\left(\log(|t|+2)+O(1)\right) +O\left(\frac{x^{-2}}{|t|+2}\right),\label{3.11}\end{split}\end{equation}
where
\begin{equation} a_n(x) = \min\left( \Big(\frac{n}{x}\Big)^{\frac{1}{2}}, \Big(\frac{x}{n}\Big)^{\frac{3}{2}}\right).\label{3.12}\end{equation} 
\end{proposition}

\noindent \textit{Proof.} This proposition is proved by using an explicit formula from Landau's 1909 Handbuch \cite{La}, which states that (unconditionally) 
for $x>1$,  $x\neq p^m$,
\begin{equation}  \sum_{n\le x}\frac{\Lambda(n)}{n^s} = -\frac{\zeta '}{\zeta}(s) + \frac{x^{1-s}}{1-s} - \sum_\rho \frac{x^{\rho-s}}{\rho-s} +\sum_{n=1}^\infty \frac{x^{-2n-s}}{2n+s}\label{3.13}\end{equation}
provided $s\neq 1$, $s\neq \rho$, $s\neq -2n$.\footnote{ If $s=0$ we get \eqref{3.1}.  Landau used \eqref{3.13} to prove Riemann's original explicit formula for $\pi(x)$.} Rewriting \eqref{3.13}, we have
\begin{equation}  \sum_\rho \frac{x^{\rho}}{\rho-s} = -x^s\Big(\frac{\zeta '}{\zeta}(s)+ \sum_{n\le x}\frac{\Lambda(n)}{n^s} \Big)  + \frac{x }{1-s}  +\sum_{n=1}^\infty \frac{x^{-2n}}{2n+s}.\label{3.14}\end{equation}
This equation holds independently of RH, but assuming RH we have $\rho = \frac{1}{2}+i\gamma$.  Letting  $s= \frac{3}{2}+it$ and using \eqref{2.2}, the above equation simplifies to read 
\[ -x^{\frac{1}{2}}\sum_{\gamma}\frac{x^{i\gamma }}{1+i(t-\gamma)} =  x^{it}\sum_{n> x}\frac{\Lambda(n)a_n(x)}{n^{it}}+ \frac{x}{-\frac{1}{2} - it}+\sum_{n=1}^\infty \frac{x^{-2n}}{2n+\frac{3}{2}+it}.\]
On the other hand, if  $s= -\frac{1}{2}+it$ in \eqref{3.14} we have
\[ \begin{split} x^{\frac{1}{2}}\sum_{\gamma}\frac{x^{i\gamma }}{1-i(t-\gamma)} &= - x^{it}\sum_{n\le x}\frac{\Lambda(n)a_n(x)}{n^{it}} - x^{-\frac{1}{2}+it}\frac{\zeta '}{\zeta}\left(-\frac{1}{2}+it\right)\\& \qquad + \frac{x}{\frac{3}{2} - it}+\sum_{n=1}^\infty \frac{x^{-2n}}{2n-\frac{1}{2}+it}.\end{split}\]
Subtracting the latter from the former and using 
\[ \frac{\zeta '}{\zeta}\left(-\frac{1}{2}+it\right) = - \log( |t|+2) + O(1) ,\]
which follows easily from the functional equation,  we obtain the proposition.  By continuity the values $x=1,p^m$ no longer need to be excluded. 
The role of RH in Proposition 1 is notational.  Recently, a new notation has emerged which is very convenient. We write the complex zeros of the zeta function as
$\rho = \frac{1}{2} +i \gamma,\ \gamma \in  \mathbb{C}$,
so that $\gamma$ is complex when that zero is off the $\frac{1}{2}$-line. Thus,
the RH becomes the statement that  $\gamma$ is real. 
With this notation we see that the proof is unchanged, and Proposition 1 holds  unconditionally.  We will not make any further use of this notation since the size  of the terms in our sums over zeros become important and the RH is often needed. 

\section{Montgomery's theorem}
We first examine Montgomery's explicit formula heuristically and see what each term means.  The weight in the sum over zeros concentrates the sum to zeros in a short bounded interval around $t$, and therefore behaves similarly to 
\[  \sum_{ t <\gamma \le t+1 }x^{i\gamma}. \]
By \eqref{2.18}, if this sum  is substantially smaller than $ \log t$ then we will have detected cancelation from $x^{i\gamma}$. If $x=1$ or is close to 1 no cancelation can occur, and this is reflected by the term $x^{-1/2}\log(|t|+2)$ in \eqref{3.11}.  The sum over primes is concentrated around $x$, and therefore behaves similarly to 
\[  \sum_{\frac{1}{2}x<n\le 2x} \frac{\Lambda(n)}{n^{it}}. \]
The expected value of the original sum over primes is obtained by the PNT  and equals the remaining term 
\[ \frac{2 x^{1-it}}{\left(\frac{1}{2}+it\right)\left(\frac{3}{2}-it\right)}.\]

How does one extract information from \eqref{3.11}? Montgomery was interested in studying the distribution of the differences of pairs of zeros, and for this it is clear one needs to square the absolute value of the sum over zeros. It would be nice to be able to obtain this distribution in an interval of length one around $t$, but the pointwise  dependence on $t$ in the Dirichlet sum over primes is intractable. To circumvent this problem, we also integrate with respect to $t$ to obtain our distribution in a longer range. To this end we consider 
\[ \int_0^T \bigg|\sum_\gamma \frac{x^{i\gamma }}{1+(t-\gamma)^2}\bigg|^2 dt. \]
Since the weight in the sum will be small when $|t-\gamma|$ is large, which is the case over most of the integration range unless $0<\gamma \le T$, we may restrict the sum to this range with a small error. With the sum restricted to the zeros $0<\gamma \le T$, we may extend the integration range to $(-\infty ,\infty)$ with a small error. Using \eqref{2.18}, Montgomery showed
\[\int_0^T \bigg|\sum_\gamma \frac{x^{i\gamma }}{1+(t-\gamma)^2}\bigg|^2 dt = \int_{-\infty}^\infty \bigg|\sum_{0<\gamma \le T} \frac{x^{i\gamma }}{1+(t-\gamma)^2}\bigg|^2 dt +O\left((\log T)^3\right).\]
Multiplying out the integral on the right-hand side, we find
\[ \int_{-\infty}^\infty \bigg|\sum_{0<\gamma \le T} \frac{x^{i\gamma }}{1+(t-\gamma)^2}\bigg|^2 dt 
=\frac{\pi}{2}\sum_{0<\gamma ,\gamma '\le T}x^{i(\gamma-\gamma ')}w(\gamma - \gamma '), \]
where the weight
\begin{equation}  w(u) = \frac{4}{4+u^2} \label{4.1}\end{equation}
is obtained on evaluating the integral either by residues, convolution, or otherwise.  
We thus define for $x>0$ 
\begin{equation} F(x,T) = \sum_{0<\gamma ,\gamma '\le T}x^{i(\gamma-\gamma ')}w(\gamma - \gamma ')=\frac{2}{\pi} \int_{-\infty}^\infty  \bigg|\sum_{0<\gamma \le T} \frac{x^{i\gamma }}{1+(t-\gamma)^2}\bigg|^2 dt 
. \label{4.2}\end{equation}
Then
\begin{equation}  F(x,T) \ge 0, \qquad F(x,T) = F\left(\frac{1}{x},T\right),\label{4.3} \end{equation}
and 
\begin{equation} F(x,T) = \frac{2}{\pi} \int_0^T \bigg|\sum_\gamma \frac{x^{i\gamma }}{1+(t-\gamma)^2}\bigg|^2 dt +O\left((\log T)^3\right).\label{4.4} \end{equation}
The next step is to use Proposition 1 to evaluate $F(x,T)$. Denoting \eqref{3.11} by
\[ L(x,t) = R(x,t),\]
we have just shown that
\begin{equation}\int_0^T\big|L(x,t)\big|^2 \, dt = 2\pi x F(x,T) + O(x(\log T)^3). \label{4.5}\end{equation}
For $R(x,T)$, we compute the mean-square of each term.  For the Dirichlet series, we use a standard mean value theorem of Montgomery and Vaughan \cite{Mo2}, which states that 
\begin{equation} \int_0^T\bigg|\sum_{n=1}^\infty\frac{a_n}{n^{it}}\bigg|^2\, dt = \sum_{n=1}^\infty |a_n|^2\big(T+O(n)\big). \label{4.6} \end{equation}
Hence 
\[ \begin{split}\int_0^T\bigg|\sum_{n=1}^\infty\frac{\Lambda(n) a_n(x)}{n^{it}}\bigg|^2 \, dt &= \sum_{n=1}^\infty \big|\Lambda(n) a_n(x)\big|^2(T +O(n))\\& = xT(\log x +O(1)) + O\left(x^2\log x\right), \end{split}\]
 by Stieltjes integration and the PNT (with remainder). 
The remaining terms are elementary:
\[ \int_0^T\bigg|\frac{2 x^{1-it}}{(\frac{1}{2}+it)(\frac{3}{2}-it)}\bigg|^2\, dt \ll x^2,\]
\[ \int_0^T\Big|x^{-\frac{1}{2}}\left(\log(|t|+2)+O(1)\right)\Big|^2\, dt = \frac{T}{x}\Big( (\log T)^2 +O(\log T)\Big) ,\]
and
\[ \int_0^T\bigg|\frac{x^{-2}}{|t|+2}\bigg|^2\, dt \ll x^{-4}.\]
We thus have two main terms, the Dirichlet series term for $ (\log T)^{3/2}\le x \le o(T)$  and the term $\log (|t|+2)$ which dominates for $1\le x \le (\log T)^{3/4}$. In the intermediate range all terms are  $o(xT\log T)$.
By the Cauchy-Schwarz inequality, the largest term among these provides the main term in an asymptotic formula. Therefore,
\[ \int_0^T |R(x,t)|^2\, dt =  xT(\log x + o(\log T)) + O(x^2\log x) + \frac{T}{x}(\log T)^2 (1+o(1)) , \]
and we conclude that
\begin{equation}  F(x,T) = \frac{T}{2\pi}\log x + o(T\log T) + O(x\log x) + \frac{T}{2 \pi x^2}(\log T)^2 (1+o(1)).\label{4.7}\end{equation}
Following Montgomery, we set
\begin{equation} x=T^\alpha \label{4.8}\end{equation}
and normalize by defining
\begin{equation}  F(\alpha) = F(\alpha ,T) = \left(\frac{T}{2\pi}\log T\right)^{-1}\sum_{0<\gamma , \gamma ' \le T} T^{i\alpha (\gamma -\gamma ')} w(\gamma -\gamma ') .\label{4.9}\end{equation}
Thus we have arrived at Montgomery's theorem.
\begin{theorem} Assume the Riemann Hypothesis. Then $F(\alpha)$ is real, even, and non-negative. Further, uniformly for $0\le \alpha \le 1-\epsilon$, we have
\begin{equation} F(\alpha) = \alpha + o(1) + (1+o(1)) T^{-2\alpha}\log T .\label{4.10}\end{equation}
\end{theorem}
The error term $O(x\log x)$  in \eqref{4.7} can be improved to $O(x)$ by using a sieve bound for prime twins \cite{GoMo}, which shows the theorem  holds for 
\begin{equation} 0\le \alpha \le 1.\label{4.11}\end{equation} 
 A detailed analysis of the above proof has recently been done by Tsz Ho Chan, with all second order terms obtained. 

\section{Application to simple zeros and small gaps between zeros}
The function $F(\alpha)$ is useful for evaluating sums over differences of zeros. Let $r(u) \in L^1$, and define the Fourier transform by
\begin{equation} \hat{r}(\alpha) = \int_{-\infty}^\infty r(u) e(\alpha u) \, du,  \qquad e(u)= e^{2\pi i u}. \label{5.1}\end{equation}
If $\hat{r}(\alpha) \in L^1$, we have almost everywhere 
\begin{equation} r(u) = \int_{-\infty}^\infty \hat{r}(\alpha) e(-u\alpha ) \, d\alpha. \label{5.2}\end{equation}
On multiply \eqref{4.9} by $\hat{r}(\alpha)$ and integrating, we obtain
\begin{equation} \sum_{0<\gamma,\gamma '\le T}r\left((\gamma -\gamma ')\frac{\log T}{2\pi}\right)w(\gamma -\gamma ') = \frac{T}{2\pi}\log T \int_{-\infty}^\infty \hat{r}(\alpha)F(\alpha) \, d\alpha . \label{5.3} \end{equation}
Using Theorem 1, we can  evaluate the right-hand side provided $\hat{r}(\alpha)$ has support in $[-1,1]$. Thus, we can evaluate sums over differences of zeros on the class of functions whose Fourier transforms are supported in $[-1,1]$. Using the Fourier pair  
\begin{equation}  k(u)= \left(\frac{\sin \pi \lambda u}{\pi \lambda u}\right)^2 , \qquad \hat{k}(\alpha)= \frac{1}{\lambda}\max\left(1-\frac{|\alpha|}{\lambda} ,0\right) \qquad (\lambda >0)\label{5.4}\end{equation}
we have for $0<\lambda \le 1$
\begin{equation}\begin{split} \sum_{0<\gamma,\gamma '\le T}&\left(\frac{\sin(\frac{\lambda}{2}(\gamma -\gamma\ ')\log T)}{\frac{\lambda}{2}(\gamma -\gamma ')\log T}\right)^2w(\gamma -\gamma ') \\&\qquad =\left( \frac{1}{\lambda}\int_{-\lambda}^\lambda \Big(1- \frac{|\alpha|}{\lambda}\Big)F(\alpha) \, d\alpha \right) \frac{T}{2\pi}\log T \\& 
\qquad \sim  \left(\frac{2}{\lambda}\int_0^\lambda \Big(1- \frac{\alpha}{\lambda}\Big)\Big(\alpha +T^{-2\alpha}\log T\Big) \, d\alpha \right) \frac{T}{2\pi}\log T \\& \qquad \sim \left( \frac{1}{\lambda}+\frac{\lambda}{3}\right)\frac{T}{2\pi}\log T .\label{5.5}\end{split}\end{equation}
This result has an important application to simple zeros of $\zeta(s)$.
\begin{theorem} Assume the Riemann Hypothesis. At least two thirds of the zeros of the Riemann zeta-function are simple in the sense that as $T\to \infty$
\begin{equation} N_s(T) := \sum_{\substack{0<\gamma \le T\\ \rho\ \mathrm{simple}}}1 \ge \left(\frac{2}{3} - o(1)\right)N(T). \label{5.6} \end{equation}
\end{theorem}
\textit{Proof.} The sum in \eqref{5.5} over pairs of zeros counts distinct zeros weighted by their multiplicity. Thus a double pole gets counted 4 times, a triple zero 9 times, etc. Denoting the multiplicity of $\rho$ by $m_\rho$, we have
\[ \begin{split} \sum_{0<\gamma \le T}m_\rho &= \sum_{\substack{0<\gamma, \gamma '\le T\\ \gamma = \gamma' }} 1 \\& \le \sum_{0<\gamma,\gamma '\le T}\left(\frac{\sin\big(\frac{\lambda}{2}(\gamma -\gamma\ '\big)\log T)}{\frac{\lambda}{2}(\gamma -\gamma ')\log T}\right)^2w(\gamma -\gamma ') \\& \le (1+o(1))\left( \frac{1}{\lambda}+\frac{\lambda}{3}\right)\frac{T}{2\pi}\log T. \end{split}\]
Choosing $\lambda =1$, we have
\begin{equation} \sum_{0<\gamma \le T}m_\rho  \le \left(\frac{4}{3} + \epsilon\right)\frac{T}{2\pi}\log T. \label{5.7}\end{equation}
But
\[\sum_{\substack{0<\gamma \le T\\ \rho\ \mathrm{simple}}}1 \ge \sum_{0<\gamma \le T}(2-m_\rho ),\]
and applying \eqref{2.17} completes the proof. 

It is possible to make very small improvements in the value $\frac{2}{3}$ in Theorem 2. It would be a major advance to be able to prove that almost all the zeros are simple, even on RH.  Conrey, Ghosh, and Gonek \cite{CGGsz} have proved using a different method that assuming RH and the Generalized Lindel\"of Hypothesis, 
\[ N_s(T) \ge \left(\frac{19}{27}-\epsilon\right) N(T). \]

Montgomery also proved that there are gaps between zeros closer than the average. He used the transform pair \eqref{5.4} with their roles reversed to obtain
\[ \liminf_{n\to \infty}\ (\gamma_{n+1}-\gamma_n) \frac{\log \gamma_n}{2\pi} \le 0.669\ldots  .\]
Consider the Fourier pair 
\begin{equation}  h(u)= \left(\frac{\sin \pi  u}{\pi  u}\right)^2 \left(\frac{1}{1-u^2}\right), \qquad \hat{h}(\alpha)= \max\left(1-|u|+\frac{\sin 2\pi |u|}{2\pi} ,0\right) ,\label{5.8}\end{equation}
where $h(u)$ is the Selberg minorant of the characteristic function of the interval  $[-1,1]$ in the class of functions with Fourier transforms with support in $[-1,1]$. We prove
\begin{theorem} Assume the Riemann Hypothesis. We have
\begin{equation} \liminf_{n\to \infty} \left(\gamma_{n+1}-\gamma_n\right) \frac{\log \gamma_n}{2\pi} \le 0.6072\ldots . \label{5.9} \end{equation}
\end{theorem}
\emph{Proof.}  Take $r(u)= h\left(\frac{u}{\lambda}\right)$. Then $r(u)$ is a minorant of the characteristic function of the interval $[-\lambda , \lambda]$. Thus
\[ \begin{split}\sum_{0< \gamma \le T}m_\rho + 2 \sum_{0< \gamma -\gamma'\le\frac{ 2 \pi \lambda}{\log T}}1 &\ge \sum_{0<\gamma,\gamma ' \le T} h\left((\gamma -\gamma ' ) \frac{\log T}{2\pi\lambda}\right)w(\gamma -\gamma ') \\& 
=\left( \frac{T}{2\pi}\log T\right)\int_{-\frac{1}{\lambda}}^{\frac{1}{\lambda}}\lambda\hat{h}(\lambda \alpha)F(\alpha)\, d\alpha .\end{split}\]
Assume $\lambda <1$. Since the integrand is positive we obtain a lower bound by decreasing the integration range to $[-1,1]$. We can assume
\[ \sum_{0< \gamma \le T}m_\rho \sim \frac{T}{2\pi}\log T ,\]
since otherwise we would have infinitely many multiple zeros and the theorem holds for this reason. Thus 
\[\sum_{0< \gamma -\gamma'\le\frac{ 2 \pi \lambda}{\log T}}1 \ge \left(\frac{1}{2}-\epsilon\right) \frac{T}{2\pi}\log T \left(\lambda -1 +2\lambda\int_0^1\alpha h(\lambda \alpha)\, d\alpha \right).\]
By an easy numerical calculation, we find that the right-hand side is positive for $\lambda > 0.6072\ldots$, which proves the result.

By a different method (on RH), Montgomery and Odlyzko \cite{M-O} improved on this result and obtained the upper bound $0.5179$. Conrey, Ghosh, and Gonek \cite{CGGzg} later replaced this by $0.5172$. 

\section{Montgomery's Conjectures}
What if $\alpha >1$?  It is not difficult to see from the proof of Montgomery's theorem that for $x\ge T$ 
\begin{equation} F(x,T)  = \frac{1}{2\pi x}\int_0^T\Big|\sum_{n=1}^\infty \frac{\Lambda(n) a_n(x)}{n^{it}} - \frac{2 x^{1-it}}{\left(\frac{1}{2}+it\right)\left(\frac{3}{2}-it\right)}\Big|^2 \, dt +o(T\log T). \label{6.1}\end{equation}
We saw that the diagonal terms in the sum contribute 
$\frac{T}{2\pi}\log x $,
while the expected value term contributes  $cx$. On the other hand, we have the trivial bound
\begin{equation} F(x,T) \le F(0,T) \sim \frac{T}{2\pi}\log^2T, \label{6.2}\end{equation}
where the last relation follows from Theorem 1 (or unconditionally from \eqref{2.14}). Thus $F(x,T)$ never gets as large as $x$ for $x\gg T(\log T)^2$, and therefore the off-diagonal terms in the sum over primes must almost perfectly cancel the expected value term.

Montgomery proceeded by multiplying out the integrand in \eqref{6.1} and integrating term by term. For the off-diagonal terms, one needs to assume the Hardy-Littlewood $k$-tuple conjecture \cite{HL} for 2-tuples (or prime pairs) with a strong error term. This conjecture states that for $0<k\le N$
\begin{equation}\sum_{n\le N}\Lambda(n)\Lambda(n+k) = \mathfrak{S} (k) N +O\left(N^{\frac{1}{2}+\epsilon}\right), \label{6.3}\end{equation}
where 
\begin{equation}
 \mathfrak{S}(k) = \left\{ \begin{array}{ll}
      {\displaystyle 2C_2\prod_{\stackrel{\scr p \vert k}{\scr p>2}}\left(\frac{p-1}{p-2}\right),} &
        \mbox{if $k$ is  even, 
	$n\neq 0$;}\\
0,   &\mbox{if $k$ is  odd;}\\
\end{array}
\right.\label{6.4} \end{equation}
and
\begin{equation}
 C_2 = \prod_{p>2}\left( 1 - \frac{1}{(p-1)^2}\right).\label{6.5}\end{equation}
Montgomery stated that this conjecture \lq \lq would allow us to carry out our program" for $x\le T\le x^{2-\epsilon}$ and obtain
\[  F(x,T) \sim \frac{T}{2\pi}\log T. \]
Further, there is no reason to expect any change in behavior for bounded $ \alpha\ge 2$. On this basis Montgomery made the following conjecture.
 
\medskip
\noindent \textbf{Strong Pair Correlation Conjecture (SPC).} For any fixed bounded  $M$,  \begin{equation} F(\alpha) = 1+o(1),  \quad \textrm{ for} \ 1\le \alpha \le M. \label{6.6}\end{equation}
 A question left unanswered by \eqref{6.6} is the rate at which the function $M=M(T)$ tends to infinity.

With regard to Montgomery's heuristics for making SPC, the argument that \eqref{6.3} implies SPC in the range $1\le \alpha \le 2-\epsilon$ was carried out by Bolanz in  a 1987 Diplomarbeit (in 131 pages).\footnote{This thesis only proves the result in the range $x\le T\le x^{\frac{3}{2}-\epsilon}$, but Bolanz extended  the result to the wider range (written communication).}  At the cost of slightly weaker but acceptable error terms, one can greatly simplify Bolanz's proof by smoothing \eqref{6.1}  (see \cite{GG}).  In section 9, we will  see that one can go further by never multiplying out the integrand in \eqref{6.1}. 

With  SPC and Theorem 1 we can now evaluate almost any sum over differences of zeros. In particular, Montgomery was lead to make the following now famous conjecture.\footnote{ The SPC conjecture doesn't explicitly say anything about pair correlation, and was often not distinguished from the PCC. It is also sometimes called Montgomery's $F(\alpha)$ conjecture.}

\medskip
\noindent \textbf{Pair Correlation Conjecture (PCC).} For any fixed $\beta >0$, 
\begin{equation}  N(T,\beta) := \left(\frac{T}{2\pi}\log T\right)^{-1} \sum_{\substack{0<\gamma,\gamma '\le T\\ 0<\gamma '-\gamma \le \frac{2\pi \beta}{\log T}}} 1\sim \int_0^\beta 1 - \left(\frac{\sin \pi u}{\pi u}\right)^2 \, du .\label{6.7}\end{equation}
The density here for the number of pairs of zeros within $\beta$ of the average spacing between zeros is where the connection with random matrix theorem first occured.

One can now replace Theorems 2 and 3 with completely satisfactory results.
From the PCC we immediately see that the following conjecture is true. 

\medskip
\noindent\textbf{Small Gaps Conjecture (SGC).} We have
\begin{equation}\liminf_{n\to \infty}\  (\gamma_{n+1} - \gamma_n )\frac{\log \gamma_n}{2 \pi} =0. \label{6.8} \end{equation}
We also have

\medskip
\noindent\textbf{Simple Zeros Conjecture (SZC).}  We have
\begin{equation} N^*(T) := \left(\frac{T}{2\pi}\log T\right)^{-1}\sum_{0< \gamma \le T} m_\rho \sim 1.\label{6.9}\end{equation}
Technically this is a conjecture on the average multiplicity which implies almost all the zeros are simple, but there is no need to make this distinction here.
Another related conjecture that follows immediately from the PCC is that 
\begin{equation}  N(T,\beta) =o(1), \qquad \mathrm{as} \ \beta \to 0^+ ;\label{6.10}\end{equation}
this conjecture and SZC together are sometimes refereed to as the Essential Simplicity Conjecture (ESC). Of course, the PCC itself implies a stronger repulsion between zeros:
as $ \beta \to 0^+ $, 
\begin{equation}  N(T,\beta) \ll \beta^3 .\label{6.11}\end{equation}

We now prove the following result.
\begin{theorem} Assume the Riemann Hypothesis.   SPC  implies PCC and SZC.
\end{theorem}

First, we need a simple consequence of Theorem 1 to handle the range when $\alpha \ge M$.

\begin{lemma} Assume the Riemann Hypothesis. We have uniformly for any $B$, possibly depending on $T$, 
\begin{equation} \int_{B}^{B+1}F(\alpha)\, d\alpha \le 3. \label{6.12}\end{equation}
\end{lemma}

\textit{Proof.} With $B=C - \frac{1}{2}$,  we have 
\[\begin{split}  \int_{C-\frac{1}{2}}^{C+\frac{1}{2}}F(\alpha)\,d\alpha 
& \le 2\int_{C-1}^{C+1}\Big(1-|\alpha - C|\Big)F(\alpha)\,d\alpha
\\& = \frac{2}{\frac{T}{2\pi}\log T} \sum_{0<\gamma,\gamma '\le T}T^{iC(\gamma-\gamma ')}\left(\frac{\sin\left(\frac{1}{2}(\gamma -\gamma\ ')\log T\right)}{\frac{1}{2}(\gamma -\gamma ')\log T}\right)^2w(\gamma -\gamma ')\\&
\le \frac{2}{\frac{T}{2\pi}\log T} \sum_{0<\gamma,\gamma '\le T}\left(\frac{\sin\left(\frac{1}{2}(\gamma -\gamma\ ')\log T\right)}{\frac{1}{2}(\gamma -\gamma ')\log T}\right)^2w(\gamma -\gamma ')\\&
\le \frac{8}{3}+\epsilon, \end{split}\]
by \eqref{5.5}.

\textit{Proof of Theorem 4.}  For SZC,  we repeat the calculation in \eqref{5.5} but now assume $\lambda \ge 1$ and use SPC for that range to find
\begin{equation}\sum_{0<\gamma,\gamma '\le T}\left(\frac{\sin\left(\frac{\lambda}{2}(\gamma -\gamma\ ')\log T\right)}{\frac{\lambda}{2}(\gamma -\gamma ')\log T}\right)^2w(\gamma -\gamma ')   \sim \left( 1+\frac{1}{3\lambda^2}\right)\frac{T}{2\pi}\log T .\label{6.13}\end{equation}
The result now follows on letting $\lambda \to \infty$.

To prove the PCC, we use the Fejer kernel from \eqref{5.4} and apply \eqref{5.3} to get
\begin{equation} \begin{split}\int_{-\infty}^\infty F(\alpha)& \bigg(\frac{\sin(\pi \beta \alpha)}{\pi \beta \alpha}\bigg)^2 \, d\alpha\\& =\left(\frac{T}{2\pi}\log T\right)^{-1} \sum_{\substack{0<\gamma,\gamma '\le T\\ |\gamma '-\gamma |\le \frac{2\pi \beta}{\log T}}} \frac{1}{\beta}\bigg( 1 -\bigg|\frac{(\gamma -\gamma ') \log T}{2\pi \beta } \bigg|\bigg)w(\gamma-\gamma')
\\&   =\frac{1}{\beta} N^*(T) + \frac{2}{\beta^2}\int_0^\beta N(T,u)\,du +O\left(\frac{\beta(1+\beta)}{(\log T)^2}\right),\label{6.14} \end{split} \end{equation}
where the error term comes from removing the factor $w(\gamma -\gamma')$.
By SZC,  $N^*(T)\sim 1$. We now evaluate the left-hand side using Theorem 1 in the range $|\alpha|\le 1$, SPC in $1<|\alpha|\le M$,  and  Lemma 1 in $|\alpha|> M$. On letting $M\to \infty$, we have
\begin{equation} \int_0^\beta N(T,u)\, du \sim \frac{\beta^2}{2} + \int_0^\beta (u-\beta) \bigg(\frac{\sin \pi u}{\pi u}\bigg)^2 \, du .\label{6.15}\end{equation}
Since
\[  \frac{1}{h}\int_{\beta -h}^\beta N(T,u)\, du\le N(T, \beta) \le \frac{1}{h}\int_{\beta }^{\beta+h} N(T,u)\, du ,\]
we obtain the PCC on differencing \eqref{6.15}.

\section{Gallagher and Mueller's Work on Pair Correlation}
A few years after Montgomery's work, Gallagher and Mueller \cite{GaMu} proved a number of interesting results on pair correlation. Their starting point is the counting function $N(T,\beta)$ in \eqref{6.7}, but rather than assuming it satisfies the PCC they assumed 
\begin{equation}  N(T,\beta) \sim \int_0^\beta 1 -\mu(\alpha) \, d\alpha , \label{7.1}\end{equation}
uniformly for $0\le \beta_0\le \beta \le \beta_1 <\infty$, as $T\to \infty$, where $\mu$ is a real, even, continuous, $L^1$ function. Thus they assumed an asymptotic density function for pair correlation, where $\mu(\alpha)$  measures the deviation from a uniform distribution corresponding to a totally random distribution of zeros. They then proved the following result. 
\begin{theorem} With $N^*(T)$  given in \eqref{6.9}, we have 
\begin{equation} \int_{-\infty}^\infty \mu(\alpha) \, d\alpha \sim  N^*(T) .\label{7.2} \end{equation}
In particular, PCC implies  SZC. 
\end{theorem}
From this  we see that 
\[\int_{-\infty}^\infty \mu(\alpha) \, d\alpha \ge 1,\]
which shows that if the zeros of the zeta-function have an asymptotic pair correlation density, then the zeros must repulse each other somewhat. Further evidence of this was later obtained  by Gallagher \cite{Ga}.

 That the PCC implies SZC follows from
\[ \int_{-\infty}^\infty\bigg(\frac{\sin \pi \alpha}{\pi \alpha}\bigg)^2 \, d\alpha =1 .\]
The notable feature here is that this result holds unconditionally. One can obtain this result on RH by first using Theorem 1 to prove \eqref{5.5}, and then evaluating the off-diagonal terms in the sum over zeros by partial summation with $N(T,\beta)$ to determine the diagonal terms $N^*(T)$.  Gallagher and Mueller replaced Theorem 1 by a result of Fujii \cite{Fujii}  (also obtained by Selberg) on $S(T)$. Let
\begin{equation} R(T, h) = \int_0^T(S(t+h)-S(t))^2 \, dt . \label{7.3}\end{equation}
Fujii proved that
\begin{equation} R(T,h) \ll T \log(2 +h\log T) \quad \text{if} \  \frac{1}{\log T} \ll h\ll 1, \label{7.4}\end{equation}
and
\begin{equation} R(T,h) \sim \frac{T}{\pi^2}\log(h\log T) \quad \text{if} \  h\log T\to \infty, \ h\ll 1. \label{7.5}\end{equation}
\emph{Proof of Theorem 5.} We have 
\[ \begin{split} \int_0^T(N(t+h)-N(t))^2 dt &= \int_0^T\Big(\sum_{t<\gamma \le t+h}1\Big)^2 \,dt\\&
\sim \left( h N^*(T) + \frac{4\pi}{\log T}\int_0^{\frac{h\log T}{2\pi}}N(T,u)\, du\right)\left(\frac{T}{2\pi}\log T\right)  .\end{split}\]
(This is \eqref{6.14} from a different perspective.)
By \eqref{2.14}, the left-hand side is also
\[ \sim T\left(h \frac{\log T}{2\pi}\right)^2 +R(T,h) .\] 
On substituting $N(T,u)$ from \eqref{7.1} and letting $h\log T \to \infty$ and $h\to 0$ so that \eqref{7.5} applies,  the theorem follows.

Gallagher and Mueller proved that for $h=\frac{2\pi \beta}{\log T}$
\begin{equation} R(T, h) \sim T\int_{-\infty}^{\infty} \min(|\alpha|,\beta) \mu(\alpha) \, d\alpha, \label{7.6} \end{equation} 
a result essentially equivalent to PCC.

Gallagher and Mueller also studied some consequences of \eqref{7.1} for primes. In particular they proved that the error in the PNT can be improved on  assuming \eqref{7.1} and RH to
\begin{equation} \psi(x)=x+o\left(x^{\frac{1}{2}}(\log x)^2\right) ,\label{7.7}\end{equation}
and  obtained an asymptotic formula for a weighted second moment for primes in short intervals first studied by Selberg \cite{Se1}. Their proof is quite complicated, since the approach in using \eqref{7.1} requires  partial summation to evaluate sums over differences of zeros, introducing many complications to handle the \lq \lq edges" of the summation. An interesting consequence  is a form of Theorem 1 obtained for $\mu(\alpha)$. 
Assuming RH and also SZC, Gallagher and Mueller proved 
\begin{equation} \hat{\mu}(\alpha) = 1 - |\alpha|, \qquad |\alpha|\le 1.\label{7.8}\end{equation}
The PCC density agrees with this,  and  has 
$\hat{\mu}(\alpha) =0$ elsewhere.

Related to this, there is  an alternative form of the PCC which has been found useful when generalizing to higher correlations. Starting from \eqref{5.3}, and supposing $\hat{r}(\alpha)$ has support in $[-1,1]$, we have by Theorem 1 on RH that
\begin{equation} \begin{split} \int_{-\infty}^\infty  \hat{r}(\alpha) F(\alpha)\, d\alpha & =  r(0)+\int_{-\infty}^\infty  \hat{r}(\alpha)\left( F(\alpha)-1\right)\, d\alpha \\& 
\sim  r(0) + \hat{r}(0)  - \int_{-1}^1(1-|\alpha|)\hat{r}(\alpha)\, d\alpha  \\&
= r(0) + \int_{-\infty}^\infty r(\alpha) \left(1- \Big(\frac{\sin \pi \alpha}{\pi \alpha}\Big)^2\right)\, d\alpha, \label{6.16}\end{split}\end{equation}
by Plancherel's formula. The second line is still true if  $F(\alpha) \sim 1$ for $|\alpha| \ge 1$ even if $\hat{r}(\alpha)$ does not have support in $[-1,1]$. It is no accident that the PCC density occurs in the integrand. The conclusion is that assuming RH, Theorem 1 implies that
\begin{equation} \sum_{0<\gamma,\gamma '\le T}r\Big((\gamma -\gamma ')\frac{\log T}{2\pi}\Big)w(\gamma -\gamma ') \sim r(0) + \int_{-\infty}^\infty r(\alpha) \left(1- \Big(\frac{\sin \pi \alpha}{\pi \alpha}\Big)^2\right)\, d\alpha, \label{6.17}\end{equation}
for all $r$ with $\hat{r}$ having support in $[-1,1]$. Moreover, PCC is equivalent to the conjecture that \eqref{6.17} holds for all test functions $r$ in some dense subset of $L^1$. Here, the factor $w(\gamma-\gamma')$ may be removed, if desired. 

\section{Heath-Brown's Results on Primes}
In \cite{HB} Heath-Brown proved a number of results on  primes using Montgomery's $F(\alpha)$ function. By \eqref{6.2}, the trivial bound for $F(\alpha)$ is
\begin{equation} F(\alpha) \le (1+o(1))\log T. \label{8.1}\end{equation}
Heath-Brown showed that any improvement in the order of magnitude of this bound would have important implications for primes. First, he proved that the improvement in the error in the PNT \eqref{7.7}  also holds if one assumes RH and $F(\alpha) = o(\log T)$ uniformly for $1\le \alpha\le M$, for any bounded $M$. It should be pointed out that further improvements in the error depend not only on the size of $F(\alpha)$ but also on the growth of $M(T)$. 

Heath-Brown next proved a number of results on gaps between primes, which take their strongest form if we assume 
\begin{equation} F(\alpha) \ll 1, \label{8.2}\end{equation}
for various ranges of $\alpha$.
With regard to Cram\'er's bound \eqref{3.9}, he proved that, assuming RH and \eqref{8.2}, for $\alpha$ in any small interval around $\alpha=2$
\begin{equation} p_{n+1}-p_n \ll \sqrt{p_n \log p_n} .\label{8.3}\end{equation}
Assuming $F(\alpha) \sim 1$ in this range one can improve  \eqref{8.3} on RH to little oh \cite{HBG}\footnote{The unusual order of the listed authors was due to a typo in the manuscript.}. (This also follows from a result in the next section.)
Next, assuming \eqref{8.2} for $1\le \alpha \le 2+\epsilon$ and RH, 
\begin{equation} \sum_{\substack{p_n\le x\\ p_{n+1}-p_n \ge H}} (p_{n+1}-p_n)\ll \frac{x\log x}{H}. \label{8.4} \end{equation}
This becomes non-trivial as soon as $\frac{H}{\log x} \to \infty$. On integrating with respect to $H$, we obtain 
\begin{equation} \sum_{p_n\le x} (p_{n+1}-p_n)^2\ll x(\log x)^2 .\label{8.5}\end{equation}
Previously  Selberg \cite{Se1}, improving on earlier work of Cram\'er \cite{Cr}, obtained these results on RH alone with an extra $\log x$ in each bound. Finally, Heath-Brown proved on RH and $F(\alpha)\sim 1$ that in any interval around $\alpha =1$
\begin{equation} \liminf_{n\to \infty} \ \left( \frac{p_{n+1}-p_n}{\log p_n}\right)=0,\label{8.6}\end{equation}
so that there exist small gaps much smaller than the average gap between primes. This result can be made to depend on the size of the error term in the asymptotic formula for $F(\alpha)$ for $\alpha$ in a neighborhood of 1. If the error term is a logarithm smaller than the main term, then one actually gets that there are infinitely often primes a bounded distance apart. 

In the next section, we shall prove these results by following a method that is structurally different but fundamentally the same as Heath-Brown's arguments. A very useful idea of Heath-Brown is the following bound for the sum over zeros \eqref{3.6}.
\begin{theorem} For $T\ge 2$,
\begin{equation} \bigg|\sum_{0\le \gamma \le T}x^{i\gamma}\bigg| \ll \sqrt{T \max_{t\le T} F(x, t)}.\label{8.7}\end{equation}
\end{theorem} 
\noindent Note that this becomes non-trivial as soon as our bound for $F$ is non-trivial. 

\noindent \emph{Proof.}   We have 
\begin{equation} F(x,T) = \int_{-\infty}^\infty \bigg|\sum_{0<\gamma\le T}x^{i\gamma}e^{i\gamma u}\bigg|^2 e^{-2|u|}\, du . \label{8.8}\end{equation}
This is related to \eqref{4.2} by Plancherel's theorem but may be  verified directly.
By Gallagher's inequality \cite{Da}, 
\[ |f(0)|\ll \int_{-1}^1 |f(u)|\, du + \int_{-1}^1 |f'(u)|\, du \]
for $f\in C^1$.
Then with 
\[ f(u) = \bigg|\sum_{0<\gamma\le T}x^{i\gamma}e^{i\gamma u}\bigg|^2 ,\] 
we obtain
\[ \begin{split} \bigg|\sum_{0\le \gamma \le T}x^{i\gamma}\bigg|^2 &\ll \int_{-1}^1 \bigg|\sum_{0<\gamma\le T}x^{i\gamma}e^{i\gamma u}\bigg|^2 \, du \\ & \qquad +\int_{-1}^1\bigg|\sum_{0<\gamma\le T}x^{i\gamma}e^{i\gamma u}\bigg| \bigg|\frac{\partial}{\partial u} \sum_{0<\gamma\le T}x^{i\gamma}e^{i\gamma u}\bigg|\, du .\end{split}\]
In the first integral on the right we insert the weight $e^{-2|u|}$ and extend the limits of integration to $(-\infty, \infty)$ to see that, by \eqref{8.8},  this is bounded by $F(x,T)$. To complete the proof, the second integral is handled similarly following an application of the Cauchy-Schwarz inequality and partial summation.

\section{Equivalence between SPC and Primes} 

In \cite{GoMo} Montgomery and I proved the following equivalence between the SPC and the second moment for primes in short intervals.
\begin{theorem} Assume the Riemann Hypothesis. If  $0<B_1\le B_2\le 1$, then
\begin{equation} I(x,\delta) :=\int_1^X\left(\psi\left((1+\delta)x\right)-\psi(x)-\delta x\right)^2\, dx \sim \frac{1}{2}\delta X^2\log\frac{1}{\delta} \label{9.1}\end{equation}
holds uniformly for $X^{-B_2}\le \delta \le X^{-B_1}$
provided
\begin{equation}F(x,T) \sim \frac{T}{2\pi}\log T \label{9.2}\end{equation}
holds for 
\[X^{B_1}(\log x)^{-3}\le T\le X^{B_2}(\log x)^3.\]
  Conversely, if $1\le A_1\le A_2 <\infty$, then \eqref{9.2} holds uniformly for $T^{A_1} \le X \le T^{A_2}$ provided that \eqref{9.1} holds uniformly for
\[X^{-\frac{1}{A_1}}(\log x)^{-3}\le \delta \le X^{-\frac{1}{A_2}}(\log x)^3.\] 
\end{theorem}
In particular, one can prove on RH that SPC is equivalent to 
\begin{equation}\int_1^X(\psi(x+h)-\psi(x)-h)^2\, dx \sim h X\log\frac{X}{h}\label{9.3}\end{equation}
for $1\le h \le X^{1-\epsilon}$,\footnote{If \eqref{9.3} holds for this range of $h$  it implies  RH.}  where an argument of Saffari and Vaughan \cite{SV} is used to move from primes in the interval $(x,x+\delta x]$ to the fixed interval $(x,x+h]$. Results 
\eqref{8.3} -- \eqref{8.6} are consequences of Theorem 7. Further, it is a straightforward exercise to show that the twin prime conjecture in the form \eqref{6.3} implies \eqref{9.1} and \eqref{9.3} in the ranges $X^{-1}\le \delta \le X^{-\frac{1}{2}-\epsilon}$ and $1\le h \le X^{\frac{1}{2} - \epsilon}$ respectively, and consequently we again obtain that \eqref{6.3} implies SPC in the range $1\le \alpha \le 2-\epsilon$. Of course, the second moment for primes in short intervals \eqref{9.1} or \eqref{9.3} is a considerably weaker hypothesis and gives the full range for SPC. 

\emph{Proof of Theorem 7.}  We follow initially the analysis in \cite{Gtex}.  Let us consider again Montgomery's explicit formula \eqref{3.11} but now aim towards obtaining a sum over primes in a short interval. This is usually done by differencing values of $x$ but Montgomery showed me the following elegant approach. Let $\kappa$, $\delta$, and $T$ be related by
\begin{equation}  e^\kappa = 1 + \delta = 1+\frac{1}{T}  \label{9.4} \end{equation}
so that $\delta= \frac{1}{T}$, and define
\begin{equation} G_\kappa(t) = \left(\frac{\sin\frac{\kappa}{2}t}{\frac{\kappa}{2}t}\right)\left( \sum_{n=1}^\infty \frac{\Lambda(n) a_n(x)}{n^{it}} -\frac{2 x^{1-it}}{\left(\frac{1}{2}+it\right)\left(\frac{3}{2}-it\right)}\right) . \label{9.5} \end{equation}
The Fourier transform of $G_\kappa(t)$ is 
\[\hat{G}_\kappa(y) =  \frac{2\pi}{\kappa} \sum_{\left|y+\frac{\log n}{2\pi}\right|<\frac{\delta}{4\pi}}\Lambda(n)a_n(x) - \frac{2\pi x}{\delta}\int_{xe^{-\kappa/2}}^{xe^{\kappa/2}}a_v\left(e^{-2\pi y}\right)\frac{dv}{v}  \]
which has the desired (but weighted) sum over primes in a short interval. By Parseval's identity, we have
\begin{equation} \int_{-\infty}^\infty \Big|G_\kappa(t)\Big|^2 \, dt = \int_{-\infty}^\infty \left|\hat{G}_\kappa(y)\right|^2 \, dy. \label{9.6}\end{equation}
Using \eqref{3.11} to express $G_\kappa(t)$ in terms of a sum over zeros with the remaining terms estimated as error terms and  simplifying we find, assuming RH, 
\begin{equation}\begin{split}
\int_0^\infty\bigg(\sum_{y<n\le y+\frac{ y}{T} }&\Lambda(n) a_n(x) - x \int_{xe^{-\kappa}}^{x}a_v(y)\, \frac{dv}{v} \bigg)^2\, \frac{dy}{y} \\&= \frac{4x\kappa^2}{\pi}\int_0^\infty \left(\frac{\sin\frac{\kappa}{2}t}{\frac{\kappa}{2}t}\right)^2\bigg|\sum_\gamma \frac{x^{i\gamma}}{1+(t-\gamma)^2}\bigg|^2\, dt +O\left(\frac{(\log T)^2}{T}\right). \label{9.7} \end{split}\end{equation}
We abbreviate this equation as
\begin{equation} L(x,T) =R(x,T).\label{9.8}\end{equation}

To prove Heath-Brown's results \eqref{8.3} and \eqref{8.4}  from the last section,  it is easy to see that, taking $T=\frac{3x}{H}$,  
\[ L(x,T) \gg \frac{x}{T^2}  \sum_{\substack{\frac{x}{2}<p_n\le x\\ p_{n+1}-p_n \ge H}} (p_{n+1}-p_n) \]
and, assuming \eqref{8.2},
\[ R(x,T) \ll \frac{x}{T}\log T.\]
Equation \eqref{8.4} follows from this, and \eqref{8.3} follows by taking only the last term in \eqref{8.4}. 

For the proof of Theorem 7 we would like to remove the weight $a_n(x)$ in $L(x,T)$ and thus obtain an expression involving $I(X,\delta)$. In view of \eqref{4.4}, since $\kappa \sim \frac{1}{T}$, $R(x,T)$ can be related to $F(x,T)$ through Abelian and Tauberian theorems. If one assumes an asymptotic formula for $I(X,\delta)$  one then obtains an asymptotic formula for $L(x,\delta)$ which gives an asymptotic formula for $R(x,\delta)$ and then a Tauberian theorem gives an asymptotic formula for $F(x,T)$.  The converse direction works similarly using an Abelian theorem. 
All the details may be found in \cite{GoMo} except how $L(x,T)$ is related to $I(X,\delta)$, since the proof there proceeds from \eqref{2.14} rather than Proposition 1. It took me a long time to figure out how to remove the weight $a_n(y)$ even though it is actually obvious. If  $\frac{y}{T}$ is small, then for $y<n\le y+ \frac{y}{T}$ it is reasonable to replace $n$ by $y$ and thus replace $a_n(x)$ with $a_y(x)$ in \eqref{9.7} with a small error. Thus the weight is removed, and one finds that 
\begin{equation} L(x,T) = 4x^3\int_x^\infty I(y,\delta)\frac{dy}{y^5} + O\left(\frac{x^2(\log T)^2}{T^3}\right).\label{9.9}\end{equation}
Since the integrand is non-negative,  if we have an asymptotic formula for $L(x,T)$ then a simple differencing argument will give an asymptotic formula for $I(x,\delta)$. The converse is immediate.  Here the error term is smaller than the main term when $T\le x \le T^{2-\epsilon}$. To obtain the full range, rather than replacing $a_n(x)$ by $a_y(x)$, we use Stieltjes integration and the PNT with  the RH error  \eqref{3.4} to evaluate the sum over primes, and together with the Cauchy Schwarz inequality we find that the error term in \eqref{9.9} can be replaced by
\[ O\left( \frac{x (\log x)^4}{T^{\frac{3}{2}}}\right) .\]
This suffices for the full range. 

\section{Selberg's theory of $S(T)$}

For more than 50 years, Selberg has been working on the distribution of values of $\log \zeta(s)$ and related functions. In the early 1940's and he made major contributions on $S(T)$ \cite{Se2,Se3}. Further results  for Dirichlet $L$-functions were obtained in \cite{Se4}. Selberg has continued to work on these problems, and while he has lectured on his results, his next published paper on this subject \cite{Se5} only appeared in 1992. In this already famous paper Selberg introduced the properties of a general class of Dirichlet series, now referred to as the \lq\lq Selberg class".   Selberg showed that his theory, originally devised for the Riemann zeta-function, carries over to the Selberg class with remarkably few changes. 
To learn more about this subject, I recommend first reading Selberg's 1992 paper.  Second, Kai-man Tsang (Selberg's only Ph.D. student) wrote a thesis \cite{Ts} in 1984 which contains full details of the proofs for some of Selberg's more recent work on $\log \zeta(s)$. Also, the two papers of A. Ghosh \cite{Gh1,Gh2}  refine  some of Selberg's work from the 1940's.  

As examples, we state two of Selberg's results proved in Tsang's thesis. Selberg has developed methods for evaluating
\[ \int_0^T F(\log \zeta(\sigma+it))\, dt,\]
for functions $F(z)$ such as sgn(Re($z$)), sgn(Im($z$)), $|\text{Re}(z)|$, and $|\text{Im}(z)|$.  Let $\chi_{\alpha,\beta}(u)$ be $1$ if $\alpha\le u\le \beta$ and zero otherwise.
Then for $\alpha < \beta$ 
\begin{equation} \int_0^T \chi_{\alpha,\beta}\left(\sqrt{\frac{\pi}{\log\log T}} S(t)\right)\, dt
= T\int_\alpha^\beta e^{-\pi u^2}\, du + O\left( \frac{ T \log\log\log T}{\sqrt{\log \log T}}\right). \label{10.1} \end{equation}
We have similar results for the real and imaginary parts of $\log \zeta(\sigma +it)$.  

For the second result, let
 $Z(T)$  denote the number of  sign changes of $S(t)$ in $[0,T]$.  Selberg proved $Z(T) \gg T(\log T)^{\frac{1}{3} - \epsilon}$ on RH in \cite{Se2}, and  unconditionally (and with an improvement  on the $\epsilon$) in \cite{Se3}.  Ghosh \cite{Gh1} improved this to  $Z(T) \gg T(\log T)^{1- \epsilon}$.\footnote{ Also obtained earlier but unpublished by Selberg}
Tsang's thesis contains the following remarkable improvements on these results. For some $c>0$, 
\begin{equation} Z(T) \gg T\log T e^{-c(\log\log\log T)^2}\label{10.2} \end{equation}
and
\begin{equation} Z(T) \ll T\log T\frac{ \log\log\log T}{\sqrt{\log\log T} }.\label{10.3} \end{equation}
If the analysis of an error term could be improved then one would obtain
\begin{equation} Z(T) \sim T \frac{\log T}{\sqrt{\pi \log\log T} }.\label{10.4} \end{equation}

 I will now describe some key ideas that went into Selberg's work on $S(t)$. The very remarkable result that Selberg proved in 1946 is that all the even moments of $S(t)$ can be computed unconditionally \cite{Se3}. He proved this on intervals $(T,T+H]$, where $T^a\le H\le T$ and $a>\frac{1}{2}$, but for simplicity we will consider the interval $[0,T]$.
\begin{theorem} For $k\ge 1$, we have
\begin{equation} \int_0^T\bigg|S(t) + \frac{1}{\pi}\sum_{p\le T^{\frac{1}{k}}}\frac{\sin(t\log p)}{\sqrt{p}}\bigg|^{2k} \, dt \ll_k T \label{10.5}\end{equation}
and
\begin{equation} \int_0^T|S(t)|^{2k}\, dt = \frac{(2k)!}{k!(2\pi)^{2k}}T(\log\log T)^{k} + O_k\left(T(\log\log T)^{k-\frac{1}{2}}\right).\label{10.6} \end{equation}\end{theorem}
This last relation is the $2k$th moment of a Gaussian. Earlier Selberg \cite{Se2}   proved \eqref{10.5}  assuming the RH,  and also \eqref{10.6} on RH but with an error term $O_k(T(\log\log T)^{k-1})$.  These results were a great advance over previous work, which had failed to even obtain an asymptotic formula for the second moment.
From \eqref{10.5}  we see that  $S(t)$ can be approximated well in $L^{2k}$ norm by the imaginary part of a short Dirichlet series. This series is short enough so that its $L^{2k}$ norm is determined by diagonal terms, and has the Gaussian property in \eqref{10.6}. Thus $S(t)$ has this property too.  

 The proof of \eqref{10.5} and \eqref{10.6} is based on an approximate formula for $S(t)$, which has its origin in Selberg's earlier paper \cite{Se1} on primes in short intervals. There, he proved on RH that for $\sigma = \frac{1}{2} +\frac{\beta}{\log T}$ and $\beta \gg 1$, 
\begin{equation}  \int_0^T\left| \frac{\zeta '}{\zeta}(\sigma +it)\right|^2\, dt \ll_\beta  T(\log T)^2 .\label{10.7} \end{equation}
 Selberg's work was ahead of its time, since we now know that replacing the bound in \eqref{10.7}  by an asymptotic formula is equivalent to the PCC \cite{GGM}. 

Selberg first found an approximate formula for $\frac{\zeta '}{\zeta}(s)$. This is not straightforward.  For $\sigma >1$, we have the Dirichlet series representation \eqref{2.2} for $\frac{\zeta '}{\zeta}(s)$. As we bring $s$ into  the critical strip the Dirichlet series fails to converge. It is a familiar fact that an appropriate partial sum of a Dirichlet series will still provide a good approximation for the analytic continuation of the series. However, on or near the critical line  we expect the poles from the zeros $\zeta(s)$ to dominate, as reflected in the partial-fraction formula, for $s\neq \rho$, $t\ge 2$,    
\begin{equation} \frac{\zeta '}{\zeta}(s) = \sum_\rho \left( \frac{1}{s-\rho}+\frac{1}{\rho}\right) +O(\log t). \label{10.8}\end{equation}
Since
\begin{equation}  S(t) = \frac{1}{\pi} \text{arg}\ \zeta\left(\frac{1}{2}+it\right) = - \frac{1}{\pi}\int_{\frac{1}{2}}^\infty \text{Im}\left( \frac{\zeta '}{\zeta}(\sigma +it)\right)\, dt ,\label{10.9} \end{equation}
it is (maybe) plausible that the Dirichlet series part of $S(t)$ will usually  dominate. A candidate for an approximate formula is \eqref{3.13} which we can rewrite as, for $x>1$, $s\neq 1$, $s\neq \rho$, $s\neq -2k$,
\begin{equation}  \frac{\zeta '}{\zeta}(s) = - \sum_{n\le x}\frac{\Lambda(n)}{n^s}+ \frac{x^{1-s}}{1-s} - \sum_\rho \frac{x^{\rho-s}}{\rho-s} +\sum_{n=1}^\infty \frac{x^{-2n-s}}{2n+s}.\label{10.10}\end{equation}
In hindsight \eqref{10.10} looks even better, because 
\[ - \frac{1}{\pi} \int_{\frac{1}{2}}^\infty \text{Im}\left(\sum_{n\le x}\frac{\Lambda(n)}{n^{\sigma+it}}\right) \, d\sigma = -\frac{1}{\pi}\sum_{n\le x}\frac{\Lambda(n) \sin(t\log n)}{\sqrt{n}\log n} \]
which gives exactly the approximation in \eqref{10.5} from the terms where $n$ is  prime.  (The prime powers will contribute an error term.) 
The problem here is that the sum over zeros does not converge absolutely, and consequently \eqref{10.10} has never been used successfully for this problem. Earlier work had smoothed this formula (or rather over-smoothed it), so that the correct approximation was lost. Selberg had the innovative idea that one only needs to smooth slightly in order to obtain  absolute convergence in the sum over zeros. 

 Let
\begin{equation}
 \Lambda_x(n)= \left\{ \begin{array}{ll}
      {\Lambda(n),} &
        \mbox{for  
	$1\le n\le x$,}\\
\Lambda(n)\frac{\log\frac{x^2}{n}}{\log n},   &\mbox{for $x\le n\le x^2$}.\\
\end{array}
\right.\label{10.11} \end{equation}
Then, for $x>1$, $s\neq 1$, $s\neq \rho$, $s\neq -2k$,
\begin{equation}  \begin{split}\frac{\zeta '}{\zeta}(s) = - \sum_{n\le x^2}\frac{\Lambda_x(n)}{n^s}+ &\frac{x^{2(1-s)}-x^{1-s}}{(1-s)^2\log x} +\frac{1}{\log x}\sum_\rho \frac{x^{\rho-s}-x^{2(\rho-s)}}{(\rho-s)^2}\\& +\frac{1}{\log x}\sum_{n=1}^\infty \frac{x^{-2n-s}-x^{-2(2n+s)}}{(2n+s)^2}.\label{10.12}\end{split}\end{equation}
This formula is much easier to prove than \eqref{10.10}.
Selberg  next argues as follows.  Assume RH,  and suppose $4\le x\le t^2$. Let
\begin{equation}  \sigma_1=\frac{1}{2}+ \frac{1}{\log x} ,\label{10.13} \end{equation}
which is at the transition from the region where the Dirichlet series dominates to the region where the zeros dominate. 
From \eqref{10.12}, we see that for $\sigma \ge \sigma_1$ and some complex number $\omega$ with $|\omega|\le 1$
\begin{equation} \frac{\zeta '}{\zeta}(\sigma +it) = - \sum_{n\le x^2}\frac{\Lambda_x(n)}{n^{\sigma+it}}+O\left(x^{\frac{1}{2}-\sigma}\right) + 2 \omega x^{\frac{1}{2}-\sigma} \sum_\gamma \frac{\sigma_1-\frac{1}{2}}{\left(\sigma_1-\frac{1}{2}\right)^2 +(t-\gamma)^2}.\label{10.14}\end{equation}
We also have, on taking the real part of \eqref{10.8},
\[ \text{Re}\frac{\zeta '}{\zeta}(\sigma +it)  = \sum_\gamma \frac{\sigma-\frac{1}{2}}{\left(\sigma-\frac{1}{2}\right)^2 +(t-\gamma)^2} +O(\log t),\]
Thus, taking the real part of \eqref{10.14} with $\sigma=\sigma_1$  gives for some $-1\le\omega '\le 1$,
\[ \begin{split} \sum_\gamma \frac{\sigma_1-\frac{1}{2}}{\left(\sigma_1-\frac{1}{2}\right)^2 +(t-\gamma)^2} +O(\log t)&= - \mathrm{Re} \sum_{n\le x^2}\frac{\Lambda_x(n)}{n^{\sigma_1+it}}\\& +O(1) +\frac{2 \omega '}{e}\sum_\gamma \frac{\sigma_1-\frac{1}{2}}{\left(\sigma_1-\frac{1}{2}\right)^2 +(t-\gamma)^2}.\end{split}\]
Since $1-\frac{2\omega '}{e} > 1 -\frac{2}{e} >\frac{1}{4}$, we conclude that
\[ \sum_\gamma \frac{\sigma_1-\frac{1}{2}}{\left(\sigma_1-\frac{1}{2}\right)^2 +(t-\gamma)^2}= O\left(\bigg|\sum_{n\le x^2}\frac{\Lambda_x(n)}{n^{\sigma_1+it}}\bigg|\right) + O(\log t).\]
Substituting this back into \eqref{10.14}, we obtain
\begin{equation} \frac{\zeta '}{\zeta}(\sigma +it) = - \sum_{n\le x^2}\frac{\Lambda_x(n)}{n^{\sigma+it}} + O\left(x^{\frac{1}{2}-\sigma}\bigg|\sum_{n\le x^2}\frac{\Lambda_x(n)}{n^{\sigma_1+it}}\bigg|\right) + O\left(x^{\frac{1}{2}-\sigma}\log t\right).\label{10.15}\end{equation}
Selberg  next substitutes \eqref{10.15} into \eqref{10.9} for the integration range $\sigma_1\le \sigma <\infty$. For  $\frac{1}{2} <\sigma \le \sigma_1$, he uses \eqref{10.8} and \eqref{10.15} to show this range only contributes to the error terms.
The conclusion is the following theorem, which is the primary tool for obtaining Theorem 8 assuming the RH.
\begin{theorem} Assume the Riemann Hypothesis. For $t\ge 2$, $4\le x\le t^2$, and $\sigma_1$ given in \eqref{10.13}, we have
\begin{equation} S(t) = - \sum_{n<x^2}\frac{\Lambda_x(n)}{n^{\sigma_1}}\frac{\sin(t\log n)}{\log n} +O\left(\frac{1}{\log x}\bigg|\sum_{n\le x^2}\frac{\Lambda_x(n)}{n^{\sigma_1+it}}\bigg|\right)+O\left(\frac{\log t}{\log x}\right). 
\label{10.16}\end{equation}
\end{theorem}

How do you remove the RH from the above analysis? I think it takes great insight to even suspect that this can be done. Selberg makes a much more subtle choice for $\sigma_1$. He defines 
\begin{equation} \sigma_{x,t} = \frac{1}{2} + 2
\max_{\rho \in \mathcal{A}} \left( \beta - \frac{1}{2}, \frac{2}{\log x}\right), \label{10.17}
\end{equation} 
where
\begin{equation} \mathcal{A} = \left\{ \rho : |t-\gamma|\le \frac{x^{3\left|\beta -\frac{1}{2}\right|}}{\log x}\right\} . \label{10.18}\end{equation}
Thus, we move towards or away from the critical line depending on how far off the line nearby zeros lie. There is also an issue of convergence, and the explicit formula \eqref{10.12} needs to be replaced by a similar formula where the sum over zeros has a factor of $(s-\rho)^3$ in the denominator. Ultimately the contribution from zeros off the $\frac{1}{2}$-line is bounded by a density estimate proved in \cite{Se3}.


\begin{thebibliography}{10}



\bibitem{BKI} E. B. Bogomolny and J. P. Keating,
\textit{Random matrix theory and the Riemann zeros I: three-
and four-point correlations},
 Nonlinearity
\textbf{8} (1995), 1115--1131.

\bibitem{BK2} E. B. Bogomolny and J. P. Keating,
\textit{Random matrix theory and the Riemann zeros II: n-point correlations},
 Nonlinearity
\textbf{9} (1996), 911--935.

\bibitem{Bo} Joachim Bolanz, \textit{Uber Die Montgomery'she Paarvermutung}, Diplomarbeit 1987, 131 pages.

\bibitem{CFKRS} J. B. Conrey, D. W. Farmer, J. P. Keating, M. O. Rubinstein, N. C. Snaith, \textit{Integral moments of L-functions}, arXiv:math.NT/0206018 , 2002, 71pp.


\bibitem{CGGzg} J. B. Conrey, A. Ghosh, and S. M. Gonek,
\textit{ A note on gaps between zeros of the  zeta-function},
 Bull. London Math. Soc. \textbf{16}
(1984), 421--424.


\bibitem{CGGsz} J. B. Conrey, A. Ghosh, and S. M. Gonek,
\textit{ Simple zeros of the Riemann  zeta-function}, Proc. London Math. Soc. (3) \textbf{76} (1998), 497--522.

\bibitem{Cr} H. Cram\'er, \textit{On the order of magnitude of the difference between consecutive prime numbers}, Acta Arithmetica, \textbf{2} (1936), 23--46.

\bibitem{Da} Harold Davenport, \textit{Multiplicative number theory.  Revised and with a preface by Hugh L. Montgomery. 3rd ed.}
Graduate Texts in Mathematics, \textbf{74 } New York, NY: Springer, 177 pp.

\bibitem{Fujii} A. Fujii,  \textit{On the zeros of Dirichlet $L$-functions, I,} Trans. A. M. S. , \textbf{196} (1974), 225--235. 

\bibitem{Fu} Akio Fujii,  \textit{On a theorem of Landau,} 
Proc. Japan Acad., Ser. A 65, No.2  (1989), 51--54. 

\bibitem{Ga} P. X. Gallagher,\textit{ Pair
correlation of zeros of the zeta function},
 J. Reine Angew. Math.
\textbf{362} (1985), 72--86.


\bibitem{GaMu} P. X. Gallagher and J. Mueller, \textit{Primes and zeros
in short intervals},
 J. Reine Angew. Math.
\textbf{ 303/304} (1978), 205--220.

\bibitem{Gh1} A. Ghosh, \textit{On Riemann's zeta function---sign changes of $S(T)$,}
Recent progress in analytic number theory, Vol. 1 (Durham, 1979),  
Academic Press, London-New York, 1981, 25--46.

\bibitem{Gh2} A. Ghosh, \textit{On the Riemann zeta function---mean value theorems and the distribution of $|S(T)|$,}
J. Number Theory \textbf{17} (1983), no. 1, 93--102.

\bibitem{Gtex} D. A. Goldston, \textit{Prime numbers and the pair correlation of zeros of the zeta-function,} Proc. of the Texas Conference on Number Theory 1982, Univ. of Texas Press. 

%\bibitem{GoCrelle} D. A. Goldston, 
%\textit{ On the pair correlation conjecture for zeros of the Riemann
%zeta-function}, J. Reine Angew. Math. \textbf{385} (1988), 24--40.



\bibitem{GG} D. A. Goldston and S. M. Gonek,
\textit{  Mean value theorems for long Dirichlet
polynomials and tails of Dirichlet series},
 Acta Arith.
\textbf{84} (1998), 155--192.

\bibitem{GGM} D. A. Goldston, S. M. Gonek, and H. L. Montgomery, 
 \textit{Mean values of the logarithmic derivative of the Riemann zeta-function with applications to primes in short intervals,}
 J. Reine Angew. Math. \textbf{537}  (2001),  105-126.
 
\bibitem{GoMo} D. A. Goldston and H. L. Montgomery
\textit{ Pair correlation of zeros 
and primes in short intervals},
 Analytic Number Theory 
and Diophantine Problems,
Birkha\"user, Boston, Mass. 1987, 183--203.

\bibitem{Gonek} S. M. Gonek, \textit{An explicit formula of Landau and its applications to the theory of the zeta-function,} Knopp, Marvin (ed.) et al., A tribute to Emil Grosswald: number theory and related analysis. Providence, RI: American Mathematical Society. Contemp. Math. 143, 395-413 (1993). 

\bibitem{HL} G. H. Hardy and J. E. Littlewood, \textit{ Some problems of `Partitio 
Numerorum': III On the expression of a number as a sum of primes}, Acta Math. 
\textbf{ 44} (1923), 1--70.

\bibitem{HB} D. R. Heath-Brown
\textit{ Gaps between primes, and the pair
correlation of   zeros of the zeta-function}
 Acta Arith.
\textbf{41} (1982), 85--99.

\bibitem{HBG} D. R. Heath-Brown and D. A. Goldston, \textit{A note on the difference between consecutive primes,} Math. Ann. \textbf{266} (1984), 317--320. 

\bibitem{He} Dennis A. Hejhal,
\textit{On the triple correlation of zeros of the zeta function}, 
 Internat. Math. Res. Notices
\textbf{ 293ff} Issue 7, (1994), 10pp (electronic).

\bibitem{Il} Georg Illies, \textit{Cram\'er functions and Guinand equations},  Acta Arith.  \textbf{105} no. 2, (2002), 103--118. 

\bibitem{In} A. E. Ingham, \textit{The distribution of prime numbers,}
Cambridge Tracts in Mathematics and Mathematical Physics, 30; Cambridge Mathematical Library, Cambridge University Press (1990), 114 pp   .

\bibitem{KS} N. Katz and P. Sarnak, \textit{Random matrices, Frobenius eigenvalues, and monodromy,} Colloquium Publications. American Mathematical Society (AMS). \textbf{45}Providence, RI, (1999) 419 pp .

\bibitem{Mo} H. L. Montgomery
\textit{ The pair correlation of zeros of the zeta function}
Proc. Sympos. Pure Math. 
\textbf{ 24} AMS, Providence, R. I., 1973, 181--193.

\bibitem{La} E. Landau, \textit{ Handbuch der Lehre von der Verteilung der Primzchlen,} Teubner, Leipzig (1909). Reprinted by Chelsea Publishing Co., New York, (1953).

\bibitem{Mo2} Hugh L. Montgomery, \textit{The analytic principle of the large sieve,}  Bull. Am. Math. Soc. \textbf{84} (1978), 547-567. 

\bibitem{M-O} H. L. Montgomery and A. Odlyzko,
\textit{ Gaps between zeros of the zeta-function},
 Topics in Classical Number Theory, Vol I, II, Budapest, 1981,
Colloquia Math. Soc. J\'anos Bolyai, \textbf{34}, North-Holland, Amsterdam-New York,
1984, 1079--1106.

%\bibitem {MS} H. L. Montgomery and K. Soundararajan, \textit{Beyond pair %correlation} to appear, 1--20.


\bibitem{O} A. M. Odlyzko,
\textit{On the distribution of spacings between zeros
of the zeta function},
 Math. Comp.
\textbf{48} (1987), 273--308.


\bibitem{RS} Ze\'ev Rudnick and Peter Sarnak,
\textit{Zeros of principal L-functions and random matrix theory.}
A celebration of John F. Nash, Jr , Duke Math. J.
\textbf{ 81} Issue 2, (1996), 269--322.

\bibitem{SV} Saffari, B. and Vaughan, R. C., \textit{On the fractional parts of $x/n$ and related sequences, II.}  Ann. Inst. Fourier (Grenoble)  \textbf{27}  (1977), no. 2,  1--30.

\bibitem{Se1} A. Selberg, \textit{On the normal density of primes in small intervals,and 
the difference between consecutive primes,} Arch. Math. Naturvid. \textbf{47} No. 6,
(1943),  87-105.

\bibitem{Se2} A. Selberg, \textit{ On the remainder in the formula for $N(T)$, the
number of zeros of $\zeta(s)$ in the strip $0<t<T$}, Avh. Norske 
Vid. Akad. Oslo. I. No.1, (1944) 27 pp.

\bibitem{Se3} A. Selberg, \textit{ Contributions to the theory of the Riemann zeta-function}, Archiv for Mathematik og Naturvidenskab B. \textbf{48} (1946), No. 5, 89--155.

\bibitem{Se4} A. Selberg, \textit{Contributions to the theory of Dirichlet's L-functions}, Skrifter utgitt av Det Norske Videnskaps-Akademi i Oslo. I. Math.-Naturv. Klasse (1946), No. 3, 1--62.

\bibitem{Se5} A. Selberg, \textit{Old and new conjectures and results about a class of Dirichlet series,} Bombieri, E. (ed.) et al., Proceedings of the Amalfi conference on analytic number theory, held at Maiori, Amalfi, Italy, Sept. 25--29, 1989. Salerno: Universitá di Salerno, (1992) 367-385 .


\bibitem{T} E. C. Titchmarsh, \textit{The theory of the Riemann zeta-function}, Second edition. Edited and with a preface by D. R. Heath-Brown. The Clarendon Press, Oxford University Press, New York, 1986.

\bibitem{Ts} Kai-man Tsang, \textit{The distribution of the values of the Riemann zeta-function,} Thesis, Princeton University, October 1984, 179pp.


\end{thebibliography}
\end{document}